\documentclass[12pt]{article}
\usepackage{latexsym}
\usepackage{amsfonts}
\usepackage{amsmath}
\usepackage{amssymb}
\usepackage{graphicx}
\usepackage{latexsym}
\usepackage{amsfonts}
\usepackage{graphicx}
\usepackage{psfrag}

\input{tcilatex}
\begin{document}

\thispagestyle{empty}

\begin{center}
{\Large \textbf{\ A note on variance bounds and location of eigenvalues }}

\vskip0.2inR. Sharma, A. Sharma and R. Saini

Department of Mathematics \& Statistics\\[0pt]
Himachal Pradesh University\\[0pt]
Shimla - 171005,\\[0pt]
India \\[0pt]
email: rajesh\_hpu\_math@yahoo.co.in
\end{center}

\vskip1.5in \noindent \textbf{Abstract. }We discuss some extensions and
refinements of the variance bounds for both real and complex numbers. The
related bounds for the eigenvalues and spread of a matrix are also derived
here.

\vskip0.5in \noindent \textbf{AMS classification. \quad } 15A42, 26C10, 60E15

\vskip0.5in \noindent \textbf{Key words and phrases}. \ Variance, complex
numbers, leptokurtic distributions, eigenvalues, trace, spread, polynomial,
span, Jung's theorem in the plane.\bigskip

\bigskip

\bigskip

\bigskip

\bigskip

\bigskip

\bigskip

\bigskip

\bigskip

\bigskip

\bigskip

\section{Introduction}

\setcounter{equation}{0}\vskip0.1inLet $z_{1},z_{2},...,z_{n}$ denote $n$
complex numbers. Their arithmetic mean is the number%
\begin{equation}
\frac{1}{n}\overset{n}{\underset{i=1}{\sum }}z_{i}=\widetilde{z}.  \tag{1.1}
\end{equation}%
In literature, the number\ \ \ \ 
\begin{equation}
\frac{1}{n}\overset{n}{\underset{i=1}{\sum }}\left\vert z_{i}-\widetilde{z}%
\right\vert ^{2}=S_{z}^{2}\text{ }  \tag{1.2}
\end{equation}%
or its equivalent expressions have been studied in several different
contexts and notations and is termed as the variance of complex numbers at
many places. For example, see Audenaert $\left[ 2\right] $, Bhatia and
Sharma $\left[ 4,5\right] ,$ Merikoski and Kumar $\left[ 13\right] $, and
Park $\left[ 17\right] $.\newline
The number%
\begin{equation}
\frac{1}{n}\overset{n}{\underset{i=1}{\sum }}\left( z_{i}-\widetilde{z}%
\right) ^{2}=S^{2}  \tag{1.3}
\end{equation}%
is also important in this context. If $z_{i}$'s are all real we denote them
by $x_{i}$'s with $a=\min x_{i}$ and $b=\max x_{i}.$ The arithmetic mean by $%
\overline{x}$ and variance by the lower case letter $s^{2}.$ In this case $%
S_{z}=\left\vert S\right\vert =S=s$ but in general $S_{z}$ rather than $%
\left\vert S\right\vert $ is more consistent with $s.$ For instance, $s=0$ $%
\left( S_{z}=0\right) $ if and only if all the $x_{i}$'s $\left( z_{i}\text{%
's}\right) $ are equal$.$ This is not the case with $\left\vert S\right\vert
;$ for example, for three distinct complex numbers $0,$ $\pm \frac{1}{2}+i%
\frac{\sqrt{3}}{2}$ we have $S=0.$ It however turns out that for some
purposes $s^{2}$ is more consistent with 
\begin{equation}
\sigma _{z}^{2}=\frac{\left\vert S^{2}\right\vert +S_{z}^{2}}{2}  \tag{1.4}
\end{equation}%
than $S_{z}^{2}.$ Note that the analogue of the Popoviciu inequality [18]%
\begin{equation}
s^{2}\leq \frac{\left( b-a\right) ^{2}}{4}  \tag{1.5}
\end{equation}%
for the complex numbers says that%
\begin{equation}
\sigma _{z}^{2}\leq \underset{i,j}{\max }\frac{\left\vert
z_{i}-z_{j}\right\vert ^{2}}{4}.  \tag{1.6}
\end{equation}%
But it is not always true that $S_{z}^{2}\leq \underset{i,j}{\max }\frac{%
\left\vert z_{i}-z_{j}\right\vert ^{2}}{4}.$ For example, for $z_{1}=-\frac{1%
}{2}+i\frac{\sqrt{3}}{2},$ $z_{2}=0$ and $z_{3}=\frac{1}{2}+i\frac{\sqrt{3}}{%
2},$ $S_{z}^{2}=\frac{1}{3}$ and $\underset{i,j}{\max }\left\vert
z_{i}-z_{j}\right\vert =1.$\newline
The corresponding inequality for $S_{z}^{2}$ is%
\begin{equation}
S_{z}^{2}\leq r_{z}^{2}\leq \underset{i,j}{\max }\frac{\left\vert
z_{i}-z_{j}\right\vert ^{2}}{3},  \tag{1.7}
\end{equation}%
where $r_{z}$ is the radius of the smallest disk containing all the numbers $%
z_{i}$'s, see $\left[ \text{4, 5}\right] .$ \newline
A classical theorem of Jung [9] says that the complex numbers $z_{i}$'s in a
plane can be contained in a closed \ disk of radius $\underset{i,j}{\max }%
\frac{\left\vert z_{i}-z_{j}\right\vert }{\sqrt{3}}.$ We thus have%
\begin{equation*}
\underset{i,j}{\max }\frac{\left\vert z_{i}-z_{j}\right\vert }{2}\leq
r_{z}\leq \underset{i,j}{\max }\frac{\left\vert z_{i}-z_{j}\right\vert }{%
\sqrt{3}}.
\end{equation*}%
In this context it is in interesting to note a case when the given complex
numbers lie on the boundary of the smallest disk containing them. We here
show that if the complex numbers lie on a circle with centre at their
arithmetic mean then this circle is the smallest circle enclosing these
points, (see Theorem 2.1 \& 3.1 below). Further, if the complex numbers are
all collinear then $\left\vert S\right\vert =\sigma _{z}=$ $S_{z}$, and
conversely, ( Theorem 2.2). A necessary and sufficient condition is given
for which the numbers $\sigma _{z}$, $S_{z}$ and $\left\vert S\right\vert $
are all equal, ( Theorem 2.2 ). We obtain a complex analogue of the
inequality, Mallows and \ Richter [11], 
\begin{equation}
s^{2}\geq \frac{r}{n-r}\left( \alpha _{r}-\overline{x}\right) ^{2}, 
\tag{1.8}
\end{equation}%
where $\alpha _{r}$ is the arithmetic mean of any subset of $r$ numbers
chosen from the real numbers $x_{1},x_{2},\ldots ,x_{n},$ (Theorem 2.3).%
\newline
On the other hand we find in literature that the inequality (1.5) and its
complementary Nagy's inequality [13], 
\begin{equation}
S^{2}\geq \frac{\left( b-a\right) ^{2}}{2n}  \tag{1.9}
\end{equation}%
also provide bounds for the spread of a complex $n\times n$ matrix $A$ when
the eigenvalues $\lambda _{i}\left( A\right) $ of $A$ are all real. The
spread of a matrix $A$ is the maximum distance between two eigenvalues of a
matrix, Spd$\left( A\right) =\lambda _{\max }\left( A\right) -\lambda _{\min
}\left( A\right) .$ We have,%
\begin{equation}
\frac{4}{n}\text{tr}B^{2}\leq \text{Spd}\left( A\right) ^{2}\leq 2\text{tr}%
B^{2},  \tag{1.10}
\end{equation}%
where $B=A-\frac{\text{tr}A}{n}I$ and tr$A$ denotes the trace of $A,$ see
[6, 23].\newline
We show that the inequalities, [3, 21],%
\begin{equation}
\frac{\left( b-a\right) ^{2}}{2n}+\frac{2}{n-2}\left( \overline{x}-\frac{a+b%
}{2}\right) ^{2}\leq s^{2}\leq \left( b-\overline{x}\right) \left( \overline{%
x}-a\right) ,  \tag{1.11}
\end{equation}%
provide some further refinements of the inequalities (1.5) and (1.9) and
consequently we get better bounds for the spread of a matrix for some
special cases, ( Theorem 2.4, 2.5, 3.2). A refinement of the inequality
(1.5) is obtained for Leptokurtic and Mesokurtic distributions, (Theorem
2.6). We obtain refinements of the eigenvalue bounds in some special cases,
(Theorem 3.3, 3.4). Likewise, the bounds for the span of a polynomial are
given, (Theorem 3.5).

\section{Main Results}

\setcounter{equation}{0}\vskip0.1in\textbf{Theorem 2.1. }If the complex
numbers $z_{i}$'s lie on a circle in the complex plane with centre $%
\widetilde{z}$ and radius $r_{z},$ then $r_{z}$ is the radius of the
smallest disk containing all the points $z_{i}$'s.\newline
\textbf{Proof. }For any complex number $c,$ we can write (1.2) in the form%
\begin{equation}
S_{z}^{2}=\frac{1}{n}\overset{n}{\underset{i=1}{\sum }}\left\vert z_{i}-c+c-%
\widetilde{z}\right\vert ^{2}=\frac{1}{n}\overset{n}{\underset{i=1}{\sum }}%
\left\vert z_{i}-c\right\vert ^{2}-\frac{1}{n}\overset{n}{\underset{i=1}{%
\sum }}\left\vert \widetilde{z}-c\right\vert ^{2}.  \tag{2.1}
\end{equation}%
Under the condition of the theorem, $\left\vert z_{i}-c\right\vert =r_{z}$
for all $i=1,2,\ldots ,n$ and therefore%
\begin{equation}
\frac{1}{n}\overset{n}{\underset{i=1}{\sum }}\left\vert z_{i}-c\right\vert
^{2}=r_{z}^{2}.  \tag{2.2}
\end{equation}%
Combining (2.1) and (2.2), we get that%
\begin{equation}
S_{z}^{2}+\left\vert \widetilde{z}-c\right\vert ^{2}=r_{z}^{2}.  \tag{2.3}
\end{equation}%
From the first inequality (1.7), $r_{z}\geq S_{z}.$ So the minimum value of
the $r_{z}$ is $S_{z}.$ This implies that if $r_{z}=S_{z}$ then $r_{z}$ is
the radius of the smallest disk containing the points $z_{i}$'s. For $%
\widetilde{z}=c,$ (2.3) gives $r_{z}=S_{z}.$ This proves the theorem. $%
\blacksquare $\newline
\newline
\textbf{Theorem 2.2. }Let $z_{1},z_{2},...,z_{n}$ be the points in the
finite complex plane and let $S_{z},$ $S$ and $\sigma _{z}$ be defined as in
(1.2), (1.3) and (1.4), respectively. Then, $S_{z}=\left\vert S\right\vert
=\sigma _{z}$ if and only if all the points $z_{1},z_{2},...,z_{n}$ lie on a
straight line.\newline
\textbf{Proof. }In the complex plane the convex combination of complex
numbers lie in the convex hull of these numbers. It follows that if the
points $z_{i}$'s are collinear then $\widetilde{z}$ also\ lies on the
straight line passing through $z_{i}$'s.\newline
From (1.2) - (1.4), we see that $S_{z}=\left\vert S\right\vert =\sigma _{z}$
if and only if 
\begin{equation}
\left\vert \overset{n}{\underset{i=1}{\sum }}\left( z_{i}-\widetilde{z}%
\right) ^{2}\right\vert =\overset{n}{\underset{i=1}{\sum }}\left\vert z_{i}-%
\widetilde{z}\right\vert ^{2}\text{ }.  \tag{2.4}
\end{equation}%
The equality occurs in triangle inequality%
\begin{equation*}
\left\vert \overset{n}{\underset{i=1}{\sum }}a_{i}\right\vert \leq \overset{n%
}{\underset{i=1}{\sum }}\left\vert a_{i}\right\vert
\end{equation*}%
if and only if the ratio of any two non-zero terms is positive that is $%
\frac{a_{i}}{a_{j}}>0,$ $\ i,j=1,2,...,n,$ see Ahlfors [1]. This means (2.4)
holds true if and only if the ratio of any two non zero terms in (2.4) is
positive, that is\newline
\ 
\begin{equation}
\left( \frac{z_{i}-\widetilde{z}}{z_{j}-\widetilde{z}}\right) ^{2}>0. 
\tag{2.5}
\end{equation}%
The square of a complex number $z$ is positive if and only \ if $z$ is real
\ and therefore (2.5) implies that $\frac{z_{i}-\widetilde{z}}{z_{j}-%
\widetilde{z}}$ is real. Also, $\frac{z_{i}-\widetilde{z}}{z_{j}-\widetilde{z%
}}$ is real if and only if $z_{i}$ lies on the straight line passing through 
$z_{j}$ and $\widetilde{z}.$ If $z_{k}-\widetilde{z}=0$ for some $k$ then $%
z_{k}=\widetilde{z}$ and so $z_{k}$ lies on the straight line passing
through $z_{j}$ and $\widetilde{z}.$ $\blacksquare $\newline
\newline
We need following lemma to extend the inequality (1.8) for complex numbers.%
\newline
\newline
\textbf{Lemma 2.1} Let $Z_{1}=\left\{ z_{1},z_{2},...,z_{n_{1}}\right\} $
and $Z_{2}=\left\{ z_{n_{1}+1},z_{n_{1}+2},...,z_{n_{1}+n_{2}}\right\} $ be
two sets of complex numbers. Denote by $\widetilde{Z}_{i}$ and $%
S_{Z_{i}}^{2} $ the arithmetic mean and variance of $Z_{i}$'s, $i=1,2,$
respectively. Then the combined variance $S_{Z_{1}\cup Z_{2}}^{2}$ of the
set $Z_{1}\cup Z_{2}$ is given by 
\begin{equation}
S_{Z_{1}\cup Z_{2}}^{2}=\frac{n_{1}}{n_{1}+n_{2}}S_{Z_{1}}^{2}+\frac{n_{2}}{%
n_{1}+n_{2}}S_{Z_{2}}^{2}+\frac{n_{1}n_{2}}{\left( n_{1}+n_{2}\right) ^{2}}%
\left\vert \widetilde{Z}_{1}-\widetilde{Z}_{2}\right\vert ^{2}.  \tag{2.6}
\end{equation}%
\textbf{Proof.} The combined variance of the set $Z_{1}\cup Z_{2}$ of $%
n_{1}+n_{2}$ numbers can be written as%
\begin{equation}
S_{Z_{1}\cup Z_{2}}^{2}=\frac{1}{n_{1}+n_{2}}\left( \underset{j=1}{\overset{%
n_{1}}{\sum }}\left\vert z_{j}-\widetilde{a}\right\vert ^{2}+\underset{%
j=n_{1}+1}{\overset{n_{1}+n_{2}}{\sum }}\left\vert z_{j}-\widetilde{a}%
\right\vert ^{2}\right) \text{ },  \tag{2.7}
\end{equation}%
where%
\begin{equation*}
\widetilde{a}=\frac{1}{n_{1}+n_{2}}\overset{n_{1}+n_{2}}{\underset{j=1}{\sum 
}}z_{j}.
\end{equation*}%
We note that 
\begin{eqnarray*}
\left\vert z_{j}-\widetilde{a}\right\vert ^{2} &=&\left\vert z_{j}-%
\widetilde{Z}_{1}+\widetilde{Z}_{1}-\widetilde{a}\right\vert ^{2} \\
&=&\left\vert z_{j}-\widetilde{Z}_{1}\right\vert ^{2}+\left\vert \widetilde{Z%
}_{1}-\widetilde{a}\right\vert ^{2}+2\func{Re}\overline{\left( \widetilde{Z}%
_{1}-\widetilde{a}\right) }\left( z_{j}-\widetilde{Z}_{1}\right) ,
\end{eqnarray*}%
\begin{equation*}
\underset{j=1}{\overset{n_{1}}{\sum }}\left( z_{j}-\widetilde{Z}_{1}\right)
=0\text{ and }\left\vert \widetilde{Z}_{1}-\widetilde{a}\right\vert =\frac{%
n_{2}}{n_{1}+n_{2}}\left\vert \widetilde{Z}_{1}-\widetilde{Z}_{2}\right\vert 
\text{ }.
\end{equation*}%
Therefore,%
\begin{equation}
\underset{j=1}{\overset{n_{1}}{\sum }}\left\vert z_{j}-\widetilde{a}%
\right\vert ^{2}=\underset{j=1}{\overset{n_{1}}{\sum }}\left\vert z_{j}-%
\widetilde{Z}_{1}\right\vert ^{2}+\frac{n_{1}n_{2}^{2}}{\left(
n_{1}+n_{2}\right) ^{2}}\left\vert \widetilde{Z}_{1}-\widetilde{Z}%
_{2}\right\vert ^{2}.  \tag{2.8}
\end{equation}%
On using similar arguments, we have%
\begin{equation}
\underset{j=n_{1}+1}{\overset{n_{1}+n_{2}}{\sum }}\left\vert z_{j}-%
\widetilde{a}\right\vert ^{2}=\underset{j=n_{1}+1}{\overset{n_{1}+n_{2}}{%
\sum }}\left\vert z_{j}-\widetilde{Z}_{2}\right\vert ^{2}+\frac{%
n_{1}^{2}n_{2}}{\left( n_{1}+n_{2}\right) ^{2}}\left\vert \widetilde{Z}_{1}-%
\widetilde{Z}_{2}\right\vert ^{2}.  \tag{2.9}
\end{equation}%
The assertions of the theorem now follow on using (2.8) and (2.9) in (2.7). $%
\blacksquare $\newline
\newline
\textbf{Theorem 2.3. }Let $\widetilde{\gamma }_{r}$ be the arithmetic mean
of any subset of $r$ numbers chosen from the set of $n$ complex numbers $%
z_{1},z_{2},\ldots ,z_{n}$ and let $\sigma _{z}^{2}$ be defined as in (1.4).
Then the inequality 
\begin{equation}
\left\vert \widetilde{\gamma }_{r}-\widetilde{z}\right\vert ^{2}\leq \frac{%
n-r}{r}\sigma _{z}^{2}\text{ }  \tag{2.10}
\end{equation}%
holds true for $1\leq r\leq n.$\newline
\textbf{Proof. }Let $Z_{1}$ and $Z_{2}$ be the disjoint sets of $r$ and $n-r$
numbers chosen from the numbers $z_{1},z_{2},\ldots ,z_{n},$ respectively.
Denote by $S_{z\left( r\right) }^{2}$ and $S_{z\left( n-r\right) }^{2}$ the
variance of $Z_{1}$ and $Z_{2},$ respectively. We now apply Lemma 2.1 and
find that%
\begin{equation}
S_{z}^{2}=\frac{r}{n}S_{z\left( r\right) }^{2}+\frac{n-r}{n}S_{z\left(
n-r\right) }^{2}+\frac{r\left( n-r\right) }{n^{2}}\left\vert \gamma
_{r}-\gamma _{n-r}\right\vert ^{2}\text{ }.  \tag{2.11}
\end{equation}%
Further,%
\begin{equation*}
\left\vert \gamma _{r}-\gamma _{n-r}\right\vert =\left\vert \gamma _{r}-%
\frac{1}{n-r}\left( \underset{i=1}{\overset{n}{\sum }}z_{i}-\underset{i=1}{%
\overset{r}{\sum }z_{i}}\right) \right\vert =\left\vert \frac{n}{n-r}\left(
\gamma _{r}-\widetilde{z}\right) \right\vert \text{ }.
\end{equation*}%
\newline
and therefore (2.1) can be written as%
\begin{equation}
S_{z}^{2}=\frac{r}{n}S_{z\left( r\right) }^{2}+\frac{n-r}{n}S_{z\left(
n-r\right) }^{2}+\frac{r}{n-r}\left\vert \gamma _{r}-\widetilde{z}%
\right\vert ^{2}\text{ }.  \tag{2.12}
\end{equation}%
\newline
On using similar arguments, we have 
\begin{equation}
S^{2}=\frac{r}{n}S_{r}^{2}+\frac{n-r}{n}S_{n-r}^{2}+\frac{r}{n-r}\left(
\gamma _{r}-\widetilde{z}\right) ^{2}\text{ }.  \tag{2.13}
\end{equation}%
On applying triangle inequality we find from (2.13) that%
\begin{equation}
\left\vert S^{2}\right\vert \geq \frac{r}{n-r}\left\vert \gamma _{r}-%
\widetilde{z}\right\vert ^{2}-\left\vert \frac{r}{n}S_{r}^{2}+\frac{n-r}{n}%
S_{n-r}^{2}\right\vert \text{ }.  \tag{2.14}
\end{equation}%
From (2.12) and (2.14), we get that%
\begin{equation}
\left\vert S^{2}\right\vert +S_{z}^{2}\geq \frac{2r}{n-r}\left\vert \gamma
_{r}-\widetilde{z}\right\vert ^{2}+\frac{r}{n}S_{z\left( r\right) }^{2}+%
\frac{n-r}{n}S_{z\left( n-r\right) }^{2}-\left\vert \frac{r}{n}S_{r}^{2}+%
\frac{n-r}{n}S_{n-r}^{2}\right\vert \text{ }.  \tag{2.15}
\end{equation}%
Again by triangle inequality, $S_{z\left( r\right) }^{2}\geq \left\vert
S_{r}^{2}\right\vert ,$ $S_{z\left( n-r\right) }^{2}\geq \left\vert
S_{n-r}^{2}\right\vert $ and therefore%
\begin{equation}
\frac{r}{n}S_{z\left( r\right) }^{2}+\frac{n-r}{n}S_{z\left( n-r\right)
}^{2}\geq \frac{r}{n}\left\vert S_{r}^{2}\right\vert +\frac{n-r}{n}%
\left\vert S_{n-r}^{2}\right\vert \geq \left\vert \frac{r}{n}S_{r}^{2}+\frac{%
n-r}{n}S_{n-r}^{2}\right\vert .  \tag{2.16}
\end{equation}%
\newline
The inequality (2.10) now follows from (2.15) and (2.16). $\blacksquare $%
\newline
\newline
The inequality (2.10) is an extension of Mallows and Richter inequality
[11]. For $r=1,$ we obtain the generalisation of the well known Samuelson's
inequality [20],%
\begin{equation*}
\sigma _{z}^{2}\geq \frac{1}{n-1}\left\vert z_{j}-\widetilde{z}\right\vert
^{2}.
\end{equation*}%
Likewise, we can prove the following extension of Nagy's inequality [13],%
\begin{equation}
\sigma _{z}^{2}\geq \frac{1}{2n}\underset{j,k}{\max }\left\vert
z_{j}-z_{k}\right\vert ^{2},\text{ \ }j,k=1,2,\ldots ,n.  \tag{2.17}
\end{equation}%
\newline
Note that\textbf{\ }for $r=1$, $S_{1}=0$ and therefore from (2.13) on using
triangle inequality we get that%
\begin{equation*}
\left\vert S_{n-1}^{2}\right\vert \leq \frac{n}{n-1}\left\vert
S^{2}\right\vert +\frac{n}{\left( n-1\right) ^{2}}\left\vert \widetilde{z}%
-z_{j}\right\vert ^{2}.
\end{equation*}%
\newline
Similarly, from (2.12), we have%
\begin{equation*}
S_{z\left( n-1\right) }^{2}=\frac{n}{n-1}S_{z}^{2}-\frac{n}{\left(
n-1\right) ^{2}}\left\vert \widetilde{z}-z_{j}\right\vert ^{2}
\end{equation*}%
and by addition we obtain the inequality%
\begin{equation*}
\sigma _{z\left( n-1\right) }^{2}=\frac{\left\vert S_{n-r}\right\vert
^{2}+S_{z\left( n-r\right) }^{2}}{2}\leq \frac{n}{n-1}\sigma _{z}^{2}.
\end{equation*}%
\newline
It then follows inductively that the inequality%
\begin{equation*}
\sigma _{z\left( m\right) }^{2}\leq \frac{n}{m}\sigma _{z}^{2},
\end{equation*}%
\newline
holds true for $m=1,2,\ldots ,n$ and therefore for $m=2,$ we have%
\begin{equation}
\sigma _{z}^{2}\geq \frac{2}{n}\sigma _{z\left( 2\right) }^{2}=\frac{2}{n}%
\left\vert z_{i}-z_{j}\right\vert ^{2}  \tag{2.18}
\end{equation}%
for all $i,j=1,2,...,n,$ $i\neq j.$ The inequality (2.18) implies (2.17).
Also, see [24]. \newline
\newline
\textbf{Theorem 2.4.} For $0\leq a<\overline{x}\leq s,$ we have%
\begin{equation}
s^{2}+\left( \frac{s^{2}-\overline{x}^{2}}{2\overline{x}}\right) ^{2}\leq 
\frac{\left( b-a\right) ^{2}}{4}  \tag{2.19}
\end{equation}%
and\bigskip\ with $n\geq 3$%
\begin{equation}
s^{2}-\frac{2}{n-2}\left( \frac{s^{2}-\overline{x}^{2}}{2\overline{x}}%
\right) ^{2}\geq \frac{\left( b-a\right) ^{2}}{2n}.  \tag{2.20}
\end{equation}%
\newline
\textbf{Proof.} The second inequality (1.11) implies \ that%
\begin{equation*}
\overline{x}^{2}\leq \left( a+b\right) \overline{x}-ab-s^{2},
\end{equation*}%
and \ therefore for $0\leq a<\overline{x},$ we can write%
\begin{equation}
\overline{x}\leq \frac{a+b}{2}-\frac{s^{2}-\overline{x}^{2}+ab}{2\overline{x}%
}\leq \frac{a+b}{2}-\frac{s^{2}-\overline{x}^{2}}{2\overline{x}}=\alpha 
\text{ }\left( \text{say}\right) .  \tag{2.21}
\end{equation}%
It is clear that $\alpha \leq $ $\frac{a+b}{2}$ and since $f\left( x\right)
=\left( x-a\right) \left( b-x\right) $ increases in the interval $\left[ a,%
\frac{a+b}{2}\right] ,a<b,$ we find that%
\begin{equation}
\left( \overline{x}-a\right) \left( b-\overline{x}\right) \leq \left( \alpha
-a\right) \left( b-\alpha \right) =\frac{\left( b-a\right) ^{2}}{4}-\left( 
\frac{\sigma ^{2}-\overline{x}^{2}}{2\overline{x}}\right) ^{2}.  \tag{2.22}
\end{equation}%
Combining (2.22) and the second inequality (1.11); we immediately get
(2.19). \newline
Further, it follows from (2.21) that for $0<\overline{x}\leq s,$ 
\begin{equation}
\left( \frac{a+b}{2}-\overline{x}\right) ^{2}\geq \left( \frac{s^{2}-%
\overline{x}^{2}}{2\overline{x}}\right) ^{2}.  \tag{2.23}
\end{equation}%
Combining (2.23) with the first inequality (1.11); a little computation
leads to (2.20). $\blacksquare $\newline
\newline
It may be noted here that the inequality (2.20) can equivalently be written
as 
\begin{equation}
\frac{m_{2}^{\prime }}{\overline{x}}\leq b-a,  \tag{2.24}
\end{equation}%
where $m_{2}^{\prime }=s^{2}+\overline{x}^{2}.$\newline
We mention an alternative proof of (2.24). From the second inequality
(1.11), 
\begin{equation}
\frac{m_{2}^{\prime }}{\overline{x}}\leq \frac{\left( a+b\right) \overline{x}%
-ab}{\overline{x}},\text{ }\overline{x}>0.  \tag{2.25}
\end{equation}%
Also, for $0\leq a<\overline{x}\leq s,$ from the inequality (1.5), we have $%
\overline{x}\leq s\leq \frac{b-a}{2}\leq \frac{b}{2}$ and for $\overline{x}%
\leq \frac{b}{2},$%
\begin{equation}
\frac{\left( a+b\right) \overline{x}-ab}{\overline{x}}\leq b-a.  \tag{2.26}
\end{equation}%
The inequality (2.24) follows from (2.25) and (2.26).\newline
\newline
\textbf{Theorem 2.5. }For $a<0$ and $2\overline{x}\geq ns,$ we have%
\begin{equation}
s^{2}+\left( \frac{\overline{x}^{2}-\frac{n}{2}s^{2}}{2\overline{x}}\right)
^{2}\leq \frac{\left( b-a\right) ^{2}}{4}  \tag{2.27}
\end{equation}%
and with $n\geq 3,$%
\begin{equation}
s^{2}-\frac{2}{n-2}\left( \frac{\overline{x}^{2}-\frac{n}{2}s^{2}}{2%
\overline{x}}\right) ^{2}\geq \frac{\left( b-a\right) ^{2}}{2n}  \tag{2.28}
\end{equation}%
\textbf{Proof.} We write (1.8) in the form%
\begin{equation}
s^{2}\geq \frac{\left( b-\overline{x}+\overline{x}-a\right) ^{2}}{2n}=\frac{%
\left( b-\overline{x}\right) ^{2}+\left( \overline{x}-a\right) ^{2}+2\left(
b-\overline{x}\right) \left( \overline{x}-a\right) }{2n}.  \tag{2.29}
\end{equation}%
Using arithmetic geometric mean inequality,%
\begin{equation}
\left( b-\overline{x}\right) ^{2}+\left( \overline{x}-a\right) ^{2}\geq
2\left( b-\overline{x}\right) \left( \overline{x}-a\right) .  \tag{2.30}
\end{equation}%
Thus, from (2.29) and (2.30),%
\begin{equation}
s^{2}\geq \frac{2}{n}\left( b-\overline{x}\right) \left( \overline{x}%
-a\right) .  \tag{2.31}
\end{equation}%
It follows from (2.31) that 
\begin{equation*}
\overline{x}^{2}\geq \left( a+b\right) \overline{x}-\frac{n}{2}s^{2}-ab
\end{equation*}%
and consequently, for $a<0$ and $\overline{x}>0,$ we have%
\begin{equation}
\overline{x}\geq \frac{a+b}{2}+\frac{1}{2\overline{x}}\left( \overline{x}%
^{2}-\frac{n}{2}s^{2}-ab\right) \geq \frac{a+b}{2}+\frac{1}{2\overline{x}}%
\left( \overline{x}-\frac{n}{2}s^{2}\right) =\beta \text{ }\left( \text{say}%
\right) .  \tag{2.32}
\end{equation}%
It is clear that $\beta \geq \frac{a+b}{2}$ for $2\overline{x}\geq ns$ and
since $f\left( x\right) =\left( \overline{x}-a\right) \left( b-\overline{x}%
\right) $ decreases in the interval $\left[ \frac{a+b}{2},b\right] ,a<b,$ we
find that%
\begin{equation}
\left( b-\overline{x}\right) \left( \overline{x}-a\right) \leq \left( \frac{%
b-a}{2}\right) ^{2}-\frac{1}{4\overline{x}}\left( \overline{x}-\frac{n}{2}%
s^{2}\right) ^{2}.  \tag{2.33}
\end{equation}%
Combining \ (2.33) with the second inequality (1.11); we immediately get
(2.27).\newline
From (2.33), we also have%
\begin{equation}
\left( \overline{x}-\frac{a+b}{2}\right) ^{2}\geq \left( \frac{\overline{x}-%
\frac{n}{2}s^{2}}{2\overline{x}}\right) ^{2}.  \tag{2.34}
\end{equation}%
The inequality (2.28) follows from (2.34) and the first inequality (1.11). $%
\blacksquare $\newline
\newline
Sharma et al. [22] have proved that 
\begin{equation}
m_{4}+3m_{2}^{2}\leq \left( b-a\right) ^{2}\left( \overline{x}-a\right)
\left( b-\overline{x}\right) ,  \tag{2.35}
\end{equation}%
where $m_{2}=s^{2}$ and $m_{4}=\frac{1}{n}\overset{n}{\underset{i=1}{\sum }}%
\left( x_{i}-\overline{x}\right) ^{4}.$\newline
If the distribution is Leptokurtic or Mesokurtic, we have, see [10],%
\begin{equation}
\frac{m_{4}}{m_{2}^{2}}\geq 3.  \tag{2.36}
\end{equation}%
We prove a refinement of the inequality (1.5) in the following theorem.%
\newline
\newline
\textbf{Theorem 2.6. }For a Leptokurtic or Mesokurtic distribution, we have%
\begin{equation}
s^{2}\leq \left( b-a\right) \sqrt{\frac{\left( \overline{x}-a\right) \left(
b-\overline{x}\right) }{6}}\leq \frac{\left( b-a\right) ^{2}}{2\sqrt{6}}. 
\tag{2.37}
\end{equation}%
\textbf{Proof. }Under the \ assumptions of the theorem the inequalities
(2.35) and (2.36) hold true. By (2.36), $3S^{4}\leq m_{4}$ and we obtain
from (2.18) that%
\begin{equation}
6s^{4}\leq \left( b-a\right) ^{2}\left( \overline{x}-a\right) \left( b-%
\overline{x}\right) .  \tag{2.38}
\end{equation}%
This gives the first inequality (2.37). The second inequality (2.37) follows
from (2.38) on using arithmetic - geometric mean inequality, $\left( 
\overline{x}-a\right) \left( b-\overline{x}\right) \leq \frac{\left(
b-a\right) ^{2}}{4}.$ $\blacksquare $\newline
We remark that the inequality (2.37) also holds true for both discrete and
continuous distributions.

\section{Bounds for eigenvalues}

\setcounter{equation}{0}\vskip0.1in Let $\mathbb{M}\left( n\right) $ denote
the algebra of all complex $n\times n$ matrices. We assume that the
eigenvalues $\lambda _{i}\left( A\right) $ of $A=\left( a_{ij}\right) \in 
\mathbb{M}\left( n\right) $ are all real, and may define respectively their
arithmetic mean and variance to be%
\begin{equation}
\overline{\lambda }\left( A\right) =\frac{1}{n}\overset{n}{\underset{i=1}{%
\sum }}\lambda _{i}\left( A\right) =\frac{\text{tr}A}{n}  \tag{3.1}
\end{equation}%
and%
\begin{equation}
s_{\lambda }^{2}=\frac{1}{n}\overset{n}{\underset{i=1}{\sum }}\left( \lambda
_{i}\left( A\right) -\overline{\lambda }\left( A\right) \right) ^{2}=\frac{%
\text{tr}A^{2}}{n}-\left( \frac{\text{tr}A}{n}\right) ^{2}=\frac{\text{tr}%
B^{2}}{n},  \tag{3.2}
\end{equation}%
where $B=A-\frac{\text{tr}A}{n}I.$ \newline
The spread of a matrix is the greatest distance between its eigenvalues. The
notion of the spread was introduced by Mirsky [14,15] and several authors
have studied bounds for the spread of a matrix, see [6,8,13,24]. \newline
\textbf{Theorem 3.1. }If trace of a unitary matrix $U\in \mathbb{M}\left(
n\right) $ is zero then the unit circle is the smallest circle enclosing the
eigenvalues of $U,$ and greatest lower bound on the Spd$\left( U\right) $ is 
$\sqrt{3}.$\newline
\textbf{Proof. }The eigenvalues of a unitary matrix $U$ all lie on the unit
circle and by assumption of the\ theorem tr$U=0.$ So, the eigenvalues $%
\lambda _{i}\left( U\right) $'s satisfy the conditions of the Theorem 2.1
and hence \ theunit circle is the smallest circle containing $\lambda
_{i}\left( U\right) $'s. It also follows from the second inequality (1.7)
that Spd$\left( U\right) \geq \sqrt{3}.$ $\blacksquare $\newline
\textbf{Example 1. }The basic circulant matrix $C$ with first row $\left(
0,1,0\ldots ,0\right) $ is a unitary matrix and its trace is zero. By
Theorem 3.1 the unit disk is the smallest disk containing eigenvalues of $C$
and Spd$C\geq \sqrt{3}.$ Also, for $n=3$ we have Spd$C=\sqrt{3}.$\newline
The following theorem is a consequence of Theorem 2.4 and provides
refinements of the inequalities (1.10). $\blacksquare $\newline
\textbf{Theorem 3.2. }Let the eigenvalues of an element $A\in \mathbb{M}%
\left( n\right) $ be all non negative and let $0<$tr$A\leq \left( n\text{tr}%
B\right) ^{\frac{1}{2}}.$ Then%
\begin{equation}
\text{Spd}\left( A\right) \geq \frac{\text{tr}A^{2}}{\text{tr}A}  \tag{3.3}
\end{equation}%
and%
\begin{equation}
\text{Spd}\left( A\right) \leq \frac{1}{\text{tr}A}\left( 2\text{tr}%
B^{2}\left( \text{tr}A\right) ^{2}-\frac{\left( n\text{tr}B^{2}-\text{tr}%
A^{2}\right) ^{2}}{n\left( n-2\right) }\right) ^{\frac{1}{2}}.  \tag{3.4}
\end{equation}%
\textbf{Proof. }Under the condition tr$A\leq \left( n\text{tr}B\right) ^{%
\frac{1}{2}}$ , we have 
\begin{equation}
\overline{\lambda }\left( A\right) =\frac{\text{tr}A}{n}\leq \left( \frac{1}{%
n}\text{tr}B^{2}\right) ^{\frac{1}{2}}=\left( \frac{\text{tr}A^{2}}{n}%
-\left( \frac{\text{tr}A}{n}\right) ^{2}\right) ^{\frac{1}{2}}=s_{\lambda }.
\tag{3.5}
\end{equation}%
Further, the eigenvalues of $A$ are all non-negative, therefore $0<\lambda
_{\min }\leq \overline{\lambda }\leq s_{\lambda }$ and Spd$\left( A\right)
=\lambda _{\max }\left( A\right) -\lambda _{\min }\left( A\right) .$ So we
can apply Theorem 2.1, the inequalities (3.3) and (3.4) follow on using
(3.1) and (3.2) in (2.10) and (2.20), respectively. $\blacksquare $\newline
\newline
\textbf{Example 2}. Let 
\begin{equation*}
A=\left[ 
\begin{array}{cccc}
1 & 1 & 1 & 1 \\ 
1 & 4 & 1 & 1 \\ 
1 & 1 & 16 & 1 \\ 
1 & 1 & 1 & 15%
\end{array}%
\right] .
\end{equation*}%
From (1.9), $81.393\leq $Spd$\left( A\right) \leq 115.11.$The matrix $A$ is
positive definite and tr$A\leq \left( n\text{tr}B^{2}\right) ^{\frac{1}{2}}.$
So, from our bounds (3.3) and (3.4) we have better estimate $85\leq $Spd$%
\left( A\right) \leq 109.77.$ $\blacksquare $\newline
Likewise, we can obtain another refinement of the inequality (1.10) on
applying Theorem 2.5. if $\lambda _{\min }\left( A\right) <0$ and $0<2$tr$%
A\leq \left( n^{3}\text{tr}B^{2}\right) ^{\frac{1}{2}},$ then 
\begin{equation}
\text{Spd}\left( A\right) \geq \frac{1}{n\text{tr}A}\left( 16n\text{tr}%
B^{2}\left( \text{tr}A\right) ^{2}+\left( 2\left( \text{tr}A\right)
^{2}-n^{2}\text{tr}B^{2}\right) ^{2}\right) ^{\frac{1}{2}}  \tag{3.6}
\end{equation}%
and%
\begin{equation}
\text{Spd}\left( A\right) \leq \frac{1}{\text{tr}A}\left( 2\text{tr}%
B^{2}\left( \text{tr}A\right) ^{2}-\frac{\left( \left( \text{tr}A\right)
^{2}-n^{2}\text{tr}B^{2}\right) ^{2}}{n\left( n-2\right) }\right) ^{\frac{1}{%
2}}.  \tag{3.7}
\end{equation}%
Further, Wolkowicz and Styan [23] have shown that if the eigenvalues of $%
A\in \mathbb{M}\left( n\right) $ are all real and $\lambda _{1}\left(
A\right) \leq \lambda _{i}\left( A\right) \leq \lambda _{n}\left( A\right) ,$
$i=1,2,\ldots ,n,$ then%
\begin{equation}
\frac{\text{tr}A}{n}-\sqrt{\frac{n-1}{n}\text{tr}B^{2}}\leq \lambda
_{1}\left( A\right) \leq \frac{\text{tr}A}{n}-\sqrt{\frac{1}{n\left(
n-1\right) }\text{tr}B^{2}}  \tag{3.8}
\end{equation}%
and%
\begin{equation}
\frac{\text{tr}A}{n}+\sqrt{\frac{n-1}{n}\text{tr}B^{2}}\leq \lambda
_{n}\left( A\right) \leq \frac{\text{tr}A}{n}+\sqrt{\frac{1}{n\left(
n-1\right) }\text{tr}B^{2}}.  \tag{3.9}
\end{equation}%
The inequalities (3.8) and (3.9) follow respectively from the inequalities,
[7,20],%
\begin{equation}
\overline{x}-\sqrt{n-1}s\leq \underset{i}{\min }x_{i}\leq \overline{x}-\frac{%
s}{\sqrt{n-1}}  \tag{3.10}
\end{equation}%
and%
\begin{equation}
\overline{x}+\frac{s}{\sqrt{n-1}}\leq \underset{i}{\max }x_{i}\leq \overline{%
x}+\sqrt{n-1}s.  \tag{3.11}
\end{equation}%
We now discuss extensions of these inequalities for the case when any one
eigenvalue of $A$ is known as in case of stochastic and singular matrices.%
\newline
It is clear from Lemma 2.1 that if $s_{n-1}^{2}$ is the variance of $n-1$
numbers obtained by excluding a number $x_{j}$ from the real numbers $%
x_{1},x_{2},\ldots ,x_{n},$ then%
\begin{equation}
s_{n-1}^{2}=\frac{n}{n-1}s^{2}-\frac{n}{\left( n-1\right) ^{2}}\left( 
\overline{x}-x_{j}\right) ^{2}.  \tag{3.13}
\end{equation}%
\textbf{Theorem 3.3. }Let the eigenvalues of $A\in \mathbb{M}\left( n\right) 
$ be all real. Let $\nu \left( A\right) $ be an eigenvalue of $A$ and denote
the remaining eigenvalues by $\nu _{i}\left( A\right) ,$ $\nu _{1}\left(
A\right) \leq \nu _{i}\left( A\right) \leq \nu _{n-1}\left( A\right) ,$ $%
i=1,2,\ldots ,n-1.$ Then, for $n\geq 3,$%
\begin{equation}
\frac{\text{tr}A-\nu \left( A\right) }{n-1}-\sqrt{n-2}s_{\nu }\leq \nu
_{1}\left( A\right) \leq \frac{\text{tr}A-\nu \left( A\right) }{n-1}-\frac{%
s_{\nu }}{\sqrt{n-2}}  \tag{3.14}
\end{equation}%
and%
\begin{equation}
\frac{\text{tr}A-\nu \left( A\right) }{n-1}+\frac{s_{\nu }}{\sqrt{n-2}}\leq
\nu _{n}\left( A\right) \leq \frac{\text{tr}A-\nu \left( A\right) }{n-1}+%
\sqrt{n-2}s_{\nu }.  \tag{3.15}
\end{equation}%
\textbf{Proof. }The arithmetic mean $\overline{\nu }\left( A\right) $ of $n-1
$ eigenvalues $\nu _{i}\left( A\right) $ can be written as%
\begin{equation}
\overline{\nu }\left( A\right) =\frac{1}{n-1}\overset{n-1}{\underset{i=1}{%
\sum }}\nu _{i}\left( A\right) =\frac{\text{tr}A-\nu \left( A\right) }{n-1}.
\tag{3.16}
\end{equation}%
By the use of (3.13) the variance of these eigenvalues is 
\begin{eqnarray}
s_{\nu }^{2} &=&\frac{1}{n-1}\overset{n-1}{\underset{i=1}{\sum }}\left( \nu
_{i}\left( A\right) -\overline{\nu }\left( A\right) \right) ^{2}=\frac{n}{n-1%
}s_{\lambda }^{2}-\frac{n}{n-1}\left( \overline{\lambda }\left( A\right)
-\nu \left( A\right) \right) ^{2}  \notag \\
&=&\frac{\text{tr}B^{2}}{n-1}-\frac{n}{n-1}\left( \frac{\text{tr}A}{n}-\nu
\left( A\right) \right) ^{2}.  \TCItag{3.17}
\end{eqnarray}%
On applying (3.10) to $n-1$ numbers $\nu _{1}\left( A\right) ,\nu _{2}\left(
A\right) ,\ldots ,\nu _{n-1}\left( A\right) $ and \ using (3.16) and (3.17);
we immediately get (3.12). Likewise, (3.15) follows \ from (3.11). $%
\blacksquare $\newline
\newline
\textbf{Theorem 3.4. }Under the conditions of the Theorem 3.3, we have%
\begin{equation}
\underset{i,j}{\max }\left\vert \nu _{i}\left( A\right) -\nu _{j}\left(
A\right) \right\vert \leq 2\left( \text{tr}B^{2}-n\left( \frac{\text{tr}A}{n}%
-\nu \left( A\right) \right) ^{2}\right)   \tag{3.18}
\end{equation}%
and%
\begin{equation}
\underset{i,j}{\max }\left\vert \nu _{i}\left( A\right) -\nu _{j}\left(
A\right) \right\vert \geq \frac{4}{n-1}\left( \text{tr}B^{2}-n\left( \frac{%
\text{tr}A}{n}-\nu _{1}\left( A\right) \right) ^{2}\right) .  \tag{3.19}
\end{equation}%
\textbf{Proof. }On using the inequalities (1.5) and (1.9), for $n-1$ numbers 
$\nu _{1}\left( A\right) ,\nu _{2}\left( A\right) ,\ldots ,\nu _{n-1}\left(
A\right) ,$ we have%
\begin{equation}
4s_{\lambda }^{2}\leq \text{Spd}\left( A\right) ^{2}=\left( \nu _{\max
}\left( A\right) -\nu _{\min }\left( A\right) \right) ^{2}\leq 2\left(
n-1\right) s_{\lambda }^{2}.  \tag{3.20}
\end{equation}%
Inserting (3.17) in (3.20), we immediately get (3.18) and (3.19) on
simplifications. $\blacksquare $\newline
\newline
\textbf{Example 3.} Let%
\begin{equation*}
A=\left[ 
\begin{array}{cccc}
1 & 2 & 9 & 4 \\ 
2 & 10 & 0 & 4 \\ 
9 & 0 & 5 & 2 \\ 
4 & 4 & 2 & 6%
\end{array}%
\right] .
\end{equation*}%
From the inequalities (3.8) , we have$-9.0688\leq \lambda _{1}\left(
A\right) \leq .644.$ The largest eigenvalue of $A$ is $16$ as all row sum of
the positive definite matrix $A$ is $16$. From (3.14) we have better
estimate for the smallest root, $-7.521\leq \lambda _{1}\left( A\right) \leq
-2.7610.\blacksquare $\newline
\newline
We now consider polynomials with real zeros. Let $f$ be a monic polynomial%
\begin{equation}
f\left( x\right) =x^{n}+a_{1}x^{n-1}+a_{2}x^{n-2}+\ldots +a_{n}  \tag{3.21}
\end{equation}%
with only real zeros. Then the length $b-a$ of the smallest interval $\left[
a,b\right] $ containing all the zeros of $f$ is called Span of $f,$ see $%
\left[ \text{12, 19}\right] .$ Denote by $D_{n}$ the span of $f$ then%
\begin{equation}
\frac{2}{n}\sqrt{\left( n-1\right) a_{1}^{2}-2na_{2}}\leq D_{n}\leq \sqrt{2%
\frac{n-1}{n}a_{1}^{2}-4a_{2}}.  \tag{3.22}
\end{equation}%
See Corollary 6.1.4 and Theorem 6.1.6 in [19].\newline
We prove a refinement of (3.22) in the following theorem.\newline
\textbf{Theorem 3.5. }Let the zeros of the polynomial (3.21) be all
non-negative and let $na_{2}\leq \left( n-2\right) a_{1}^{2}.$ Then%
\begin{equation}
D_{n}\geq \sqrt{\frac{2a_{2}-a_{1}^{2}}{a_{1}}}  \tag{3.23}
\end{equation}%
and with $n\geq 3,$%
\begin{equation}
D_{n}\leq \sqrt{\frac{2}{n}\left( \left( n-1\right) a_{1}^{2}-2na_{2}\right)
-\frac{1}{n\left( n-2\right) }\left( \frac{2na_{2}-\left( n-2\right)
a_{1}^{2}}{a_{1}}\right) ^{2}}.  \tag{3.24}
\end{equation}%
\textbf{Proof. }Let $x_{1},x_{2},\ldots ,x_{n}$ be the roots of the
polynomial (3.21). Then, on using relation between roots and coefficient of
polynomial, we have%
\begin{equation*}
\overline{x}=\frac{1}{n}\sum x_{i}=\frac{-a_{1}}{n}
\end{equation*}%
and%
\begin{eqnarray*}
s^{2} &=&\frac{1}{n}\sum x_{i}^{2}-\overline{x}^{2}=\frac{1}{n}\left( \sum
x_{i}\right) ^{2}-\frac{2}{n}\underset{i<j}{\sum }x_{i}x_{j}-\overline{x}^{2}
\\
&=&\frac{\left( n-1\right) a_{1}^{2}-na_{2}}{n^{2}}.
\end{eqnarray*}%
The assertions of the theorem now follow on applying Theorem 2.4.

\ \ \


\begin{thebibliography}{99}
\bibitem{1} Ahlfors, L. V., \textit{Complex Analysis}, McGraw-Hill (3rd ed)
(1979).

\bibitem{2} Audenaert K.M.R., \textit{Variance bounds with an application to
norm bounds for commutators}, Linear Algebra Appl. 432 (2010), 1126-1143.\ 

\bibitem{3} Bhatia, R. and Davis, C., \textit{A better bound on the Variance}%
, Amer. Math. Month., 107 (2000), 353-357.

\bibitem{4} Bhatia, R., Sharma, R., \textit{Some inequalities for positive
linear maps}, Linear Algebra Appl. 436 (2012), 1562-1571.

\bibitem{5} Bhatia, R., Sharma, R., \textit{Positive linear maps and Spread
of matrices-II}, Linear Algebra Appl. 491 (2016), 30-40.

\bibitem{6} Brauer, A., Mewborn, A.C., \textit{The greatest distance between
two characteristic roots of a matrix}, Duke Math. J., 26 (1959), 653-661.

\bibitem{7} Brunk, H.D., \textit{Note on two papers of K.R. Nair}, J. Indian
Soc. Agricultural Statist. 11 (1959), 186-189.

\bibitem{8} Johnson , C.R., Kumar, R. and Wolkowicz, H., \textit{Lower
bounds for the spread of a matrix}, Linear Algebra Appl. 71 (1985), 161-173.

\bibitem{9} Jung, heinrich, \textit{Uber die Kleinste Kugel, die eine
raumliche Figureinschlie}$\beta $\textit{t} . J.Reine Angew Math. 123
(1901), 241-257.

\bibitem{10} Kendal, M., Stuart, A., \textit{The Advanced theory of
Statistics}, Charles Griffin \& Company Limited, Vol. 1, (1977).

\bibitem{11} Mallow, C. L., Richter, D., \textit{Inequalities of Chebyshev
type involving conditional expectations}, Ann. Math. Statist. 40 (1969),
1922-1932.

\bibitem{12} Marden, M., \textit{Geometry of Polynomials}, Amer. Math. Soc.,
(2005).

\bibitem{13} Merikoski, K. J., Kumar, R., \textit{Characterization and lower
bounds for the spread of a normal matrix}, Linear Algebra Appl. 364 (2003),
13-31.

\bibitem{14} Mirsky, L., \textit{The spread of a matrix}, Mathematika 3
(1956), 127-130.

\bibitem{15} Mirsky, L., \textit{Inequalities for normal and hermitian
matrices}, Duke Math. j. 24 (1957), 591-598.

\bibitem{16} Nagy, J. V. S., \textit{Uber algebraische Gleichungen mit
lauter reelen Wurzeln} $,$ Jahresbericht der Duetschen Mathematiker -
Vereinigung, 27 (1918), 37-43.

\bibitem{17} Park, Kun Il. \textit{Fundamentals of Probability and
Stochastic Processes with Applications to Communications}. Springer, (2018).

\bibitem{18} Popoviciu, T., \textit{Sur les equations algebriques ayant
toutes leurs racines reeles}, Mathematica (Cluj), \ 9 (1935), 129-145.

\bibitem{19} Rahman, Q.I., Schmeisser, G., \textit{Analytic Theory of
Polynomials}, Oxford University Press, (2002).

\bibitem{20} Samuelson, P. A., \textit{How deviant can you be?}, J. Amer.
Statist. Assoc., 63 (1968), 1522-1525.

\bibitem{21} Sharma, R., Shandil, R.G., Devi S, Ram, S., Kapoor, G. and
Barnett, N.S., \textit{Some bounds on the Sample Variance in terms of the
Mean and Extreme values}, \textquotedblleft Advances in Inequalities from
Probability Theory \& Statistics\textquotedblright . Edited by N.S. Barnett
\& S.S. Dragomir, 187-193, (2008).

\bibitem{22} Sharma, R., Kumar, R. and Saini, R., \textit{Bounds on spread
of matrices related to fourth central moment-II}, Pre-print.

\bibitem{23} Wolkowicz, H. and Styan, P. H., \textit{Bounds for Eigenvalues
Using Traces}, Linear Algebra and Its Appl., 29 (1980), 471-506.

\bibitem{24} Zhang, P. and Yang, H., \textit{Improvments \ in the upper
bounds for the spread of a matrix}, Math. Inequal. Appl. 18 (2015), 337-345.
\end{thebibliography}
\end{document}